\documentclass{amsart}%
\usepackage{amsfonts}
\usepackage{amsmath}
\usepackage{amssymb}
\usepackage{graphicx}%
\newtheorem{theorem}{Theorem}
\theoremstyle{plain}

\newtheorem{corollary}[theorem]{Corollary}

\newtheorem{lemma}[theorem]{Lemma}

\newtheorem{proposition}[theorem]{Proposition}

\numberwithin{equation}{section} \numberwithin{theorem}{section}

\begin{document}
\title{Analytic regularity of a free boundary problem}
\date{28 August 2005}
\author{Huiqiang Jiang}
\address{School of Mathematics, University of Minnesota, 127 Vincent Hall, 206 Church
St. S.E., Minneapolis, MN 55455} \email{hqjiang@math.umn.edu}

\begin{abstract}
In this paper, we consider a free boundary problem with volume
constraint. We show that positive minimizer is locally Lipschitz and
the free boundary is analytic away from a singular set with
Hausdorff dimension at most $n-8$.
\end{abstract}

\maketitle

\section{Introduction}

Let $\Omega$ be a bounded open domain in $\mathbb{R}^{n}$, $n\geq2$. We use
$\mathcal{M}_{\Omega}$ to denote the collection of all pairs of $\left(
A,u\right)  $ such that $A\subset\Omega$ is a set of finite perimeter and
$u\in H^{1}\left(  \Omega\right)  $ satisfies
\[
u\left(  x\right)  =0\text{ a.e. }x\in A.
\]

We consider the energy functional
\begin{equation}
E_{\Omega}\left(  A,u\right)  =\int_{\Omega}\left\vert \triangledown
u\right\vert ^{2}+P_{\Omega}\left(  A\right)  ,\label{Equation energy}%
\end{equation}
defined on $\mathcal{M}_{\Omega}$, where $P_{\Omega}\left(  A\right)  $
denotes the perimeter of $A$ inside $\Omega$ in the sense of De Giorgi, i.e.,%
\[
P_{\Omega}\left(  A\right)  =\mathcal{H}^{n-1}\left(  \partial^{\ast}%
A\cap\Omega\right)  ,
\]
where $\partial^{\ast}A$ is the reduced boundary of $A$ and $\mathcal{H}%
^{n-1}$ is the $\left(  n-1\right)  $-dimensional Hausdorff measure.

A pair $\left(  A,u\right)  \in\mathcal{M}_{\Omega}$ is said to be a local
minimizer of $\left(  \ref{Equation energy}\right)  $ in its volume class if
for any $\left(  \tilde{A},\tilde{u}\right)  \in\mathcal{M}_{\Omega}$ such
that $\left(  \tilde{A},\tilde{u}\right)  $ agrees with $\left(  A,u\right)  $
away from a compact set and satisfies the volume constraint $\left\vert
\tilde{A}\right\vert =\left\vert A\right\vert $, we have%
\[
E_{\Omega}\left(  A,u\right)  \leq E_{\Omega}\left(  \tilde{A},\tilde
{u}\right)  .
\]
And we say $\left(  A,u\right)  $ is a nonnegative local minimizer of $\left(
\ref{Equation energy}\right)  $ in its volume class if in addition $u$ is nonnegative.

This free boundary problem is a special case of what was considered in
\cite{2004Jiang_Lin} where $u$ maps $\Omega$ to $\mathbb{R}^{p}$, $p\geq1$ and
$u\left(  A\right)  \subset\Sigma$ where $\Sigma$ is a smooth submanifold in
$\mathbb{R}^{p}$. Hence, all the results in \cite{2004Jiang_Lin} hold.
Especially, $\partial^{\ast}A$ satisfies the so called mass ratio lower bound,
i.e., given $K\subset\subset\Omega$, there exists $r_{K}>0$ and $\lambda
_{K}>0$, such that for any $x\in K\cap\overline{\partial^{\ast}A}$, and for
any $r<r_{K}$,%
\[
\mathcal{H}^{n-1}\left(  \partial^{\ast}A\cap B_{r}\left(  x\right)  \right)
\geq\lambda_{K}r^{n-1}.
\]
Let%
\[
A^{\ast}=\left\{  x\in\Omega\backslash\overline{\partial^{\ast}A}:\left\vert
B_{r}\left(  x\right)  \cap A\right\vert =\left\vert B_{r}\left(  x\right)
\right\vert \text{ for some }r>0\right\}  ,
\]
a consequence of mass ratio lower bound of $\partial^{\ast}A$ is that the
symmetric difference $A^{\ast}\triangle A$ has $\mathcal{H}^{n}$-measure zero,
i.e., the open set $A^{\ast}$ is equivalent to $A$ as a set of finite
perimeter. So we can assume $A=A^{\ast}$. Let $\partial A$ be the topological
boundary of the open set $A$, mass ratio lower bound of $\partial^{\ast}A$ now
implies%
\[
\mathcal{H}^{n-1}\left(  \left(  \partial A\backslash\partial^{\ast}A\right)
\cap\Omega\right)  =0.
\]

Now we can state our main result:

\begin{theorem}
\label{Theorem main result}Let $\left(  A,u\right)  $ be a nonnegative local
minimizer of $\left(  \ref{Equation energy}\right)  $ in its volume class.
Then $u$ is locally Lipschitz in $\Omega$ and $\partial^{\ast}A\cap\Omega$ is
an analytic hypersurface with $\mathcal{H}^{s}\left(  \left(  \partial
A\backslash\partial^{\ast}A\right)  \cap\Omega\right)  =0$ for any $s>\left(
n-8\right)  $.
\end{theorem}

If we drop the nonnegative assumption in the above theorem, it is proved in
\cite{2004Jiang_Lin} that $u\in C^{\frac{1}{2}}\left(  \Omega\right)  $, and
the proof of Theorem \ref{Theorem main result} implies that $\partial^{\ast
}A\cap\Omega$ is an analytic hypersurface away from the set where $u$ changes
sign. On the other hand, when the space dimension is two, the full regularity
in the sign changing case has recently been obtained in a joint work with C.
Larsen\cite{0000Jiang_Larsen} where a totally different blow up argument was used.

We also remark that the volume constraint is not essential for our regularity
results. Locally, volume constraint is of higher order than $P_{\Omega}\left(
A\right)  $ so it disappears after blowing up. If we drop the volume
constraint or if we instead add a volume term $c\left\vert A\right\vert $ in
the energy, the results of Theorem \ref{Theorem main result} still hold.

A related problem was considered in
\cite{2001Athanasopoulos_Caffarelli_Kenig_Salsa} by I. Athanasopoulos, L. A.
Caffarelli, C. Kenig and S. Salsa. Given $g\in H^{1}\left(  \Omega\right)  $,
they were interested in the minimizer of $E_{\Omega}\left(  A,u\right)  $
where $u\in H^{1}\left(  \Omega\right)  $ satisfies the boundary condition
$u-g\in H_{0}^{1}\left(  \Omega\right)  $ and $A$ is the set such that
$u\geq0$ in $A$ and $u\leq0$ in $\Omega\backslash A$. The free boundary
problem we are considering is quite different from theirs. Nonetheless, their
techniques in proving the Lipschitz continuity of $u$ still work for our
nonnegative local minimizer.

The paper is organized in the following way: First, we collect some results
proved in \cite{2004Jiang_Lin} and deduce the positive density property of the
free boundary. In section $\ref{Section Lipschitz}$, we prove the Lipschitz
continuity of $u$ following the arguments in
\cite{2001Athanasopoulos_Caffarelli_Kenig_Salsa}. Finally, we show
$\partial^{\ast}A\cap\Omega$ is analytic away from a singular set with
Hausdorff dimension at most $n-8$ by deriving the Euler-Lagrange equation of
the free boundary using domain variation.

\section{Preliminaries}

To prove the regularity of free boundary using a variational approach, we need
to construct good candidates to compare with. The volume constrain adds
difficulty to such construction, luckily, we can ignore the volume constraint
as long as we are willing to pay some penalty. More precisely, let $\left(
A,u\right)  $ be a local minimizer of $\left(  \ref{Equation energy}\right)  $
in its volume class, we have

\begin{lemma}
\label{Lemma deformation}There exists $r_{0}>0$, such that for any
$$x\in \Omega,\,\, r<\min\left\{  r_{0},\operatorname*{dist}\left(
x,\partial \Omega\right)  \right\}  $$ and for any pair $\left(
A_{1},u_{1}\right)  $ which agrees with $\left(  A,u\right)  $ away
from $B_{r}\left(  x\right)  $, there exists $\left(
A_{2},u_{2}\right)  $ which agrees with $\left( A_{1},u_{1}\right) $
in $B_{r}\left(  x\right)  $ and agrees with $\left( A,u\right)  $
away from a precompact subset of $\Omega$ such that $\left\vert
A_{2}\right\vert =\left\vert A\right\vert $, $u_{2}\left(  x\right)
=0$ a.e.
$x\in A_{2}$ and%
\begin{equation}
E_{\Omega}\left(  A_{2},u_{2}\right)  \leq E_{\Omega}\left(  A_{1}%
,u_{1}\right)  +C\left\vert \left\vert A\right\vert -\left\vert A_{1}%
\right\vert \right\vert \label{Equation deformation error}%
\end{equation}
for some positive constant $C$ independent of $x$ and $r$.
\end{lemma}

Before going to the proof of Lemma \ref{Lemma deformation}, let's first recall
a deformation lemma. We write any point $x\in\mathbb{R}^{n}$ as $x=\left(
x^{\prime},x_{n}\right)  \in\mathbb{R}^{n-1}\times\mathbb{R}$. Let $\xi$ be a
smooth function on $\left[  0,1\right)  $ such that $0\leq\xi\left(  r\right)
\leq1$, $\xi\left(  r\right)  \equiv1$ if $0\leq r\leq\frac{1}{4}$,
$\xi\left(  r\right)  \equiv0$ if $r\geq\frac{1}{2}$ and $\left\vert
\xi^{\prime}\left(  r\right)  \right\vert \leq8$ for any\ $r\in\left[
0,1\right)  $. For any $\varepsilon\in\mathbb{R}$, we introduce a map
$f_{\varepsilon}$ from $B_{1}(0)$ into $\mathbb{R}^{n}$ defined by%
\[
f_{\varepsilon}\left(  x\right)  =f_{\varepsilon}(x^{\prime},x_{n}%
)=(x^{\prime},x_{n}+\varepsilon\xi(|x^{\prime}|)\xi(|x_{n}|)),
\]
then we have

\begin{lemma}
\label{Lemma estimate} There exists a positive constant $\varepsilon
_{0}=\varepsilon_{0}\left(  n\right)  $ such that for any $\varepsilon
\in\left(  -\varepsilon_{0},\varepsilon_{0}\right)  $, $f_{\varepsilon}$ is a
diffeomorphism from $B_{1}(0)\subset\mathbb{R}^{n}$ into itself satisfying the
following estimates:

\begin{enumerate}
\item For any $u\in H^{1}\left(  B_{1}\left(  0\right)  \right)  $, we have%
\[
(1-c_{1}|\varepsilon|)\int_{B_{1}(0)}\left\vert \triangledown u\right\vert
^{2}\leq\int_{B_{1}(0)}\left\vert \triangledown u_{\varepsilon}\right\vert
^{2}\leq(1+c_{1}|\varepsilon|)\int_{B_{1}(0)}\left\vert \triangledown
u\right\vert ^{2}%
\]
where $u_{\varepsilon}=u\circ f_{\varepsilon}$ and $c_{1}$ is a positive
constant depending only on $n$.

\item Let $A\subset B_{1}\left(  0\right)  $ be a set of finite perimeter,
then we have%
\begin{align*}
& \left\vert \mathcal{H}^{n-1}\left(  \partial^{\ast}A_{\varepsilon}\cap
B_{1}\left(  0\right)  \right)  -\mathcal{H}^{n-1}\left(  \partial^{\ast}A\cap
B_{1}\left(  0\right)  \right)  \right\vert \\
\leq & c_{2}|\varepsilon|\mathcal{H}^{n-1}\left(  \partial^{\ast}A\cap
B_{1}\left(  0\right)  \right)  ;
\end{align*}
where $A_{\varepsilon}=f_{\varepsilon}(A)$ and $c_{2}$ is a positive constant
depending only on $n$.

\item Let $A\subset B_{1}\left(  0\right)  $ be such that the
symmetric difference $A\triangle B_{1}^{+}(0)$ satisfying%
\[
\left\vert A\triangle B_{1}^{+}(0)\right\vert \leq\delta_{0}%
\]
where $B_{1}^{+}\left(  0\right)  =\left\{  x\in B_{1}(0):x_{n}>0\right\}  $
and $\delta_{0}$ is a positive constant depending only on $n$. Then for some
positice constant $c_{3}$ depending only on $n$, we have%
\[
|A_{\varepsilon}|\leq(1-c_{3}|\varepsilon|)|A|
\]
if $\varepsilon>0$ and%
\[
|A_{\varepsilon}|\geq(1+c_{3}|\varepsilon|)|A|
\]
if $\varepsilon<0$.
\end{enumerate}
\end{lemma}

We refer the readers to \cite{2004Jiang_Lin}\cite{1999Lin_Kohn}\cite{1993Lin}
for its proof.

\begin{proof}
[Proof of Lemma \ref{Lemma deformation}]Since $A$ is a set of finite
perimeter, by a theorem of De Giorgi, $\partial^{\ast}A$ is $\left(
n-1\right)  $-rectifiable, and for every $x\in\partial^{\ast}A$, there is a
hyperplane $\Pi$ passing $x$ such that, if we denote $H^{\pm}$, the two half
spaces in $\mathbb{R}^{n}$ separated by $\Pi$, then as $r\rightarrow0^{+}$
\[
r^{-n}\left\Vert \chi_{H^{+}}-\chi_{A}\right\Vert _{L^{1}\left(  B_{r}\left(
x\right)  \right)  }\rightarrow0.
\]
After a rotation if necessary, we can always assume
$$\Pi=H_{0}=\left\{ x\in\mathbb{R}^{n}\text{: }x_{n}=0\right\}.$$
Let $r_{0}$ be sufficiently small, there exist $r_{1}>0$ and finite
number of balls $$B_{r_{k}}\left( x_{k}\right) \subset\Omega,\,\,
1\leq k\leq K$$ such that $x_{k}\in
\partial^{\ast}A$, and for any $B_{r}\left(  x\right)  \subset\Omega$,
$r<r_{0}$, there exists $k$, such that%
\begin{equation}
B_{r_{1}}\left(  x_{k}\right)  \cap B_{r}\left(  x\right)  =\emptyset
.\label{Equation empty}%
\end{equation}
After a scaling if necessary, we can assume $r_{1}=1$, and with respect to the
tangent plane $\Pi_{k}$ of $\partial^{\ast}A$ at $x_{k}$,
\[
\left\Vert \chi_{H_{k}^{+}}-\chi_{A}\right\Vert _{L^{1}\left(  B_{1}\left(
x_{k}\right)  \right)  }\leq\delta_{0},
\]
where $\delta_{0}$ is defined in Lemma \ref{Lemma estimate}. And we further
assume $r_{0}$ is so small such that
\[
r_{0}^{n}\leq c\varepsilon_{0}%
\]
for some small number $c$ depending on $n$. Now let $\left(  A_{1}%
,u_{1}\right)  $ be a pair which agrees with $\left(  A,u\right)  $ away from
$B_{r}\left(  x\right)  \subset\Omega$ with $r<r_{0}$ and $B_{1}\left(
x_{k}\right)  $ be the ball such that $\left(  \ref{Equation empty}\right)  $
holds. We define the new pair $\left(  A_{2},u_{2}\right)  $ so that it agrees
with $\left(  A_{1},u_{1}\right)  $ away from $B_{1}\left(  x_{k}\right)  $,
and we can deform inside $B_{1}\left(  x_{k}\right)  $ to meet the volume
constraint from the estimates in the third part of Lemma \ref{Lemma estimate}.
Finally, $\left(  \ref{Equation deformation error}\right)  $ follows from the
estimates in the first two parts of Lemma \ref{Lemma estimate}.
\end{proof}

Next, we recall the H\"{o}lder continuity of $u$ and the mass ratio lower
bound of $\partial^{\ast}A$ proved in \cite{2004Jiang_Lin}:

\begin{proposition}
\label{Proposition MRLB}$u\in C^{\frac{1}{2}}(\Omega)$ and $\partial^{\ast}A$
satisfies mass ratio lower bound. i.e., for any $K\subset\subset\Omega$, there
are constants $r_{K},\lambda_{K}>0$, such that for all $x\in K\cap
\overline{\partial^{\ast}A}$, $0<r<r_{K}$,
\[
\mathcal{H}^{n-1}\left(  \partial^{\ast}A\cap B_{r}\left(  x\right)  \right)
\geq\lambda_{K}r^{n-1}.
\]

\end{proposition}

Mass ratio lower bound of $\partial^{\ast}A$ implies

\begin{lemma}
\label{Lemma open}%
\[
\mathcal{H}^{n-1}\left(  \left(  \overline{\partial^{\ast}A}\backslash
\partial^{\ast}A\right)  \cap\Omega\right)  =0.
\]

\end{lemma}

\begin{proof}
Since $\partial^{\ast}A\cap\Omega$ is $\mathcal{H}^{n-1}$ measurable and
$\mathcal{H}^{n-1}\left(  \partial^{\ast}A\cap\Omega\right)  <\infty$,
standard density lemma implies that
\[
\lim_{r\rightarrow0^{+}}\frac{\mathcal{H}^{n-1}\left(  \partial^{\ast}A\cap
B_{r}\left(  x\right)  \right)  }{\omega_{n-1}r^{n-1}}=0
\]
holds for $\mathcal{H}^{n-1}$ a.e. $x\in\Omega\backslash\partial^{\ast}A$. On
the other hand, Lemma $\left(  \ref{Proposition MRLB}\right)  $ implies that
for any $x\in\overline{\partial^{\ast}A}\cap\Omega$,%
\[
\liminf_{r\rightarrow0^{+}}\frac{\mathcal{H}^{n-1}\left(  \partial^{\ast}A\cap
B_{r}\left(  x\right)  \right)  }{\omega_{n-1}r^{n-1}}>0,
\]
hence we conclude%
\[
\mathcal{H}^{n-1}\left(  \left(  \overline{\partial^{\ast}A}\backslash
\partial^{\ast}A\right)  \cap\Omega\right)  =0.
\]
\end{proof}

Define%
\begin{equation}
A^{\ast}=\left\{  x\in\Omega\backslash\overline{\partial^{\ast}A}%
:\lim_{r\rightarrow0^{+}}\frac{\left\vert B_{r}\left(  x\right)  \cap
A\right\vert }{\left\vert B_{r}\left(  x\right)  \right\vert }=1\right\}  ,
\label{Equation open set}%
\end{equation}
then $A^{\ast}$ is an open set with%
\[
\partial A^{\ast}\cap\Omega\subset\overline{\partial^{\ast}A}\cap\Omega
\]
where $\partial A^{\ast}$ is the topological boundary of $A^{\ast}$.

\begin{lemma}
The symmetric difference $A^{\ast}\triangle A$ has $\mathcal{H}^{n}$-measure
zero, hence $A^{\ast}$ and $A$ are equivalent as sets of finite perimeter.
\end{lemma}

\begin{proof}
Let $O$ be any connected open set such that $\partial^{\ast}A\cap O$ is empty,
then it is well known that either $\left\vert O\backslash A\right\vert =0$ or
$\left\vert O\cap A\right\vert =0$. Apply this observation to each component
of the open set $A^{\ast}$, we conclude $\left\vert A^{\ast}\backslash
A\right\vert =0$. Similarly,
\[
B^{\ast}=\left\{  x\in\Omega\backslash\overline{\partial^{\ast}A}%
:\lim_{r\rightarrow0^{+}}\frac{\left\vert B_{r}\left(  x\right)  \cap
A\right\vert }{\left\vert B_{r}\left(  x\right)  \right\vert }=0\right\}
\]
is also an open set such that $\partial^{\ast}A\cap B^{\ast}=\emptyset$, and
we can deduce $\left\vert B^{\ast}\cap A\right\vert =0$. From Lemma
\ref{Lemma open}, $\overline{\partial^{\ast}A}\cap\Omega$ has finite
$\mathcal{H}^{n-1}$-measure and hence zero $\mathcal{H}^{n}$-measure. Since
$\Omega$ is the disjoint union of $A^{\ast},B^{\ast}$ and $\overline{\partial
A^{\ast}}\cap\Omega$, we have%
\[
\left\vert A^{\ast}\triangle A\right\vert =\left\vert A^{\ast}\backslash
A\right\vert +\left\vert A\backslash A^{\ast}\right\vert =\left\vert A^{\ast
}\backslash A\right\vert +\left\vert B^{\ast}\cap A\right\vert +\left\vert
\overline{\partial A^{\ast}}\cap A\right\vert =0.
\]
\end{proof}

From now on, we always assume that $A$ is the open set defined by
\ref{Equation open set}. Let $\partial A$ be the topological boundary of $A$,
then it is easy to verify that%
\[
\partial A\cap\Omega=\overline{\partial^{\ast}A}\cap\Omega,
\]
hence%
\[
\mathcal{H}^{n-1}\left(  \left(  \partial A\backslash\partial^{\ast}A\right)
\cap\Omega\right)  =0.
\]

Another application of mass ratio lower bound is the following positive
density lemma which we will use in the proof of Lipschitz continuity of $u$:

\begin{lemma}
\label{Lemma positive density}For any closed set $K\subset\subset\Omega$,
there exists a constant $\lambda_{K}^{\prime}>0$, such that for any
$x\in\partial A\cap K$, and for any $r\leq\frac{1}{2}\operatorname*{dist}%
\left(  K,\partial\Omega\right)  $, we have
\[
\left\vert A\cap B_{r}\left(  x\right)  \right\vert \geq\lambda_{K}^{\prime
}r^{n}.
\]

\end{lemma}

\begin{proof}
If it is not true, then there would be a sequence $B_{r_{k}}\left(
x_{k}\right)  $, such that $x_{k}\in\partial A\cap K$, $r_{k}\leq\frac{1}%
{2}\operatorname*{dist}\left(  K,\partial\Omega\right)  $ while
\begin{equation}
\left\vert A\cap B_{r_{k}}\left(  x_{k}\right)  \right\vert \leq\frac{1}%
{k}r_{k}^{n}. \label{Equation go to zero}%
\end{equation}

First we claim
\[
\lim_{k\rightarrow\infty}r_{k}=0.
\]
Otherwise, since $\partial A\cap K$ is a compact set, extracting a subsequence
if necessary, we can assume%
\[
\lim_{k\rightarrow\infty}x_{k}=x_{0}\in\partial A\cap K
\]
and%
\[
\lim_{k\rightarrow\infty}r_{k}=r_{0}>0.
\]
From $\left(  \ref{Equation go to zero}\right)  $, we also have%
\[
\left\vert A\cap B_{r_{0}}\left(  x_{0}\right)  \right\vert =0
\]
which contradicts to $x\in\partial A$.

Next, we choose $\rho_{k}\in(\frac{r_{k}}{2},r_{k})$, such that
\[
\mathcal{H}^{n-1}(\partial^{\ast}(A\setminus B_{\rho_{k}}\left(  x_{k}\right)
)\cap\partial B_{\rho_{k}}\left(  x_{k}\right)  )\leq c(n)\frac{1}{k}%
r_{k}^{n-1}.
\]
Let $A_{k}=A\setminus B_{\rho_{k}}$, since $(A,u)$ is minimizing, when $k$ is
large , we have $r_{k}$ is small, applying Lemma \ref{Lemma deformation}, we
have,%
\[
E_{\Omega}(A,u)\leq E_{\Omega}(A_{k},u)+C\left\vert \left\vert A\right\vert
-\left\vert A_{k}\right\vert \right\vert =E_{\Omega}(A_{k},u)+C\left\vert
A\cap B_{\rho_{k}}\left(  x_{k}\right)  \right\vert ,
\]
where the last term came from the penalty for volume constraint. Hence%
\begin{align*}
&  \mathcal{H}^{n-1}\left(  \partial^{\ast}A\cap\overline{B_{\rho_{k}}\left(
x_{k}\right)  }\right) \\
\leq &  \mathcal{H}^{n-1}\left(  \partial^{\ast}(A\setminus B_{\rho_{k}%
}\left(  x_{k}\right)  )\cap\partial B_{\rho_{k}}\left(  x_{k}\right)
\right)  +C\left\vert A\cap B_{\rho_{k}}\left(  x_{k}\right)  \right\vert \\
\leq &  c(n)\frac{1}{k}r_{k}^{n-1}+C\frac{1}{k}r_{k}^{n}%
\end{align*}
which contradicts the mass ratio lower bound when $k$ is chosen sufficiently large.
\end{proof}

\section{Lipschitz continuity of $u$\label{Section Lipschitz}}

In this section, we will show that $u$ is locally Lipschitz continuous in
$\Omega$ using the approach in
\cite{2001Athanasopoulos_Caffarelli_Kenig_Salsa}. Let $\left(  A,u\right)  $
be a nonnegative local minimizer of $\left(  \ref{Equation energy}\right)  $
in its volume class. Our first step is to show that $u$ grows at most linearly
near the free boundary.

\begin{lemma}
Let $B_{4r^{\ast}}\left(  x^{\ast}\right)  \subset\Omega$ be such that
$\partial A\cap B_{r^{\ast}}\left(  x^{\ast}\right)  \neq\emptyset$, then for
some positive constant $C$,%
\[
u\left(  x\right)  \leq C\operatorname*{dist}\left(  x,\partial A\right)
\]
holds for any $x\in B_{\frac{1}{2}r^{\ast}}\left(  x^{\ast}\right)  \backslash
A$.
\end{lemma}

\begin{proof}
Let $\varphi\in C_{0}^{\infty}\left(  \mathbb{R}^{n}\right)  $ be a cutoff
function such that $0\leq\varphi\left(  x\right)  \leq1$ for any
$x\in\mathbb{R}^{n}$, $\varphi\equiv1$ in $B_{\frac{1}{2}r^{\ast}}\left(
x^{\ast}\right)  $, $\varphi\equiv0$ outside $B_{r^{\ast}}\left(  x^{\ast
}\right)  $. For any $\varepsilon\in\left(  0,\left\Vert u\right\Vert
_{L^{\infty}\left(  B_{r^{\ast}}\left(  x^{\ast}\right)  \right)  }\right)  $,
we define
\[
w=(u-\varepsilon)^{+},
\]
which is a continuous function in $B_{r^{\ast}}\left(  x^{\ast}\right)  $. Now
we consider%
\[
M=\sup_{x\in B_{r^{\ast}}\left(  x^{\ast}\right)  \backslash A}\frac{w\left(
x\right)  \varphi\left(  x\right)  }{d\left(  x\right)  }%
\]
where $d\left(  x\right)  =\operatorname*{dist}\left(  x,\partial A\right)  $.
It is easy to see that $M$ is finite and it is achieved at some point
$x_{0}\in B_{r^{\ast}}\left(  x^{\ast}\right)  $, i.e.,%
\[
Md\left(  x_{0}\right)  =w\left(  x_{0}\right)  \varphi\left(  x_{0}\right)
,
\]
and since $\partial A\cap B_{r^{\ast}}\left(  x^{\ast}\right)  \neq\emptyset$,
we have $d\left(  x_{0}\right)  =\left\vert x_{0}-y_{0}\right\vert <2r^{\ast}$
for some $y_{0}\in\partial A\cap B_{3r^{\ast}}\left(  x^{\ast}\right)  $. By a
rotation and translation if necessary, we may assume that $y_{0}=0$ and
$x_{0}=d\left(  x_{0}\right)  e_{1}$, and we also write%
\[
x=\left(  x_{1},x^{\prime}\right)  .
\]
Let%
\[
Q\left(  x-x_{0}\right)  =\frac{1}{2}\left(  x-x_{0}\right)  D^{2}\left(
w\varphi\right)  \left(  x_{0}\right)  \left(  x-x_{0}\right)  ^{\top},
\]
where $D^{2}$ is the Hessian matrix, and%
\[
\overline{Q}\left(  x^{\prime}\right)  =Q\left(  d\left(  x_{0}\right)
,x^{\prime}\right)  .
\]
Using the maximality of $\frac{w\left(  x\right)  \varphi\left(  x\right)
}{d\left(  x\right)  }$ at $x_{0}$, we can show%
\begin{equation}
\triangle_{x^{\prime}}\overline{Q}\left(  x^{\prime}\right)  \geq-\frac
{CM}{\varphi\left(  x_{0}\right)  } \label{Equation curvature lower bound}%
\end{equation}
and on the hyperplane $x_{1}=d\left(  x_{0}\right)  $,
\begin{equation}
d\left(  x\right)  \geq d\left(  x_{0}\right)  +\frac{\overline{Q}\left(
x^{\prime}\right)  }{M}+O\left(  \frac{\left\vert x^{\prime}\right\vert ^{3}%
}{M}\right)  . \label{Equation bound of free boundary}%
\end{equation}
We refer the readers to the proof of $\left(  4.3\right)  ,\left(  4.4\right)
$ in \cite{2001Athanasopoulos_Caffarelli_Kenig_Salsa} for more details. Hence,
near the origin, the free boundary $\partial A$ is below the surface
\[
S=\left\{  \left(  x_{1},x^{\prime}\right)  :x_{1}=\psi\left(  x^{\prime
}\right)  =-\frac{\overline{Q}\left(  x^{\prime}\right)  }{M}+\frac
{C\left\vert x^{\prime}\right\vert ^{3}}{M}\right\}  .
\]
Let $\kappa_{S}$ be the mean curvature of $S$, positive if convex with respect
to $e_{1}$, we have%
\[
\kappa_{S}\left(  0\right)  =-\frac{1}{n-1}\triangle\psi\left(  0\right)
\leq\frac{C}{\varphi\left(  x_{0}\right)  },
\]
hence for any $x=\left(  \psi\left(  x^{\prime}\right)  ,x^{\prime}\right)
\in S$ with $\left\vert x^{\prime}\right\vert $ small, we have%
\[
\kappa_{S}\left(  x\right)  \leq\frac{C}{\varphi\left(  x_{0}\right)
}+O\left(  \left\vert x^{\prime}\right\vert \right)  .
\]
Next, we define two families of surfaces
\[
S_{t}^{-}=\left\{  \left(  x_{1},x^{\prime}\right)  :x_{1}=\psi_{t}^{-}\left(
x^{\prime}\right)  =\psi\left(  x^{\prime}\right)  +\frac{\delta_{0}}%
{\varphi\left(  x_{0}\right)  }\left\vert x^{\prime}\right\vert ^{2}%
-t\right\}  ,
\]
and
\[
S_{t}^{+}=\left\{  \left(  x_{1},x^{\prime}\right)  :x_{1}=\psi_{t}^{+}\left(
x^{\prime}\right)  =\psi\left(  x^{\prime}\right)  +t\right\}
\]
where $t>0$, $\delta_{0}>0$ both are small. Denote by $Z_{t}$ the lens-shaped
domain between $S_{t}^{+}$ and $S_{t}^{-}$, i.e.,
\[
Z_{t}=\left\{  \psi_{t}^{-}\left(  x^{\prime}\right)  <x_{1}<\psi_{t}%
^{+}\left(  x^{\prime}\right)  \right\}  ,
\]
then $Z_{t}\subset B_{3r^{\ast}}\left(  x^{\ast}\right)  $ when $t$ is
sufficiently small. Let%
\[
V_{t}=A\cap\left\{  \psi_{t}^{-}\left(  x^{\prime}\right)  <x_{1}<\psi\left(
x^{\prime}\right)  \right\}  .
\]
We define a competing pair $\left(  A_{t},u_{t}\right)  $ such that inside
$B_{3r^{\ast}}\left(  x^{\ast}\right)  $, $A_{t}=A\backslash Z_{t}$ and%
\[
u_{t}=\left\{
\begin{array}
[c]{ccc}%
u & \text{in} & B_{3r}\backslash Z_{t},\\
v_{t} & \text{in} & Z_{t}%
\end{array}
\right.
\]
where $v_{t}$ is the harmonic extension of $u$ in $Z_{t}$, i.e., $v_{t}$ is
harmonic in $Z_{t}$ and $v_{t}=u$ on $\partial Z_{t}$. We also apply Lemma
\ref{Lemma deformation} away from $B_{3r^{\ast}}\left(  x^{\ast}\right)  $ to
keep the volume constraint which produces an extra energy of size at most
$C\left\vert V_{t}\right\vert $. Hence, since $\left(  A,u\right)  $ is a
minimizer, we have
\begin{equation}
\int_{Z_{t}}\left\vert \triangledown u\right\vert ^{2}+P\left(  A,B_{3r^{\ast
}}\left(  x^{\ast}\right)  \right)  \leq\int_{Z_{t}}\left\vert \triangledown
v_{t}\right\vert ^{2}+P\left(  A_{t},B_{3r^{\ast}}\left(  x^{\ast}\right)
\right)  +C\left\vert V_{t}\right\vert . \label{Equation energy comparison}%
\end{equation}
Next, we claim that near the origin,
\begin{equation}
u\left(  x\right)  \geq\frac{cM}{\varphi\left(  x_{0}\right)  }x_{1}+o\left(
\left\vert x\right\vert \right)  .
\label{Equation harmonic function positivity}%
\end{equation}
In fact, $u$ is positive and harmonic in $B_{d\left(  x_{0}\right)  }\left(
x_{0}\right)  $ and
\[
u\left(  x_{0}\right)  \geq\varepsilon+\frac{Md\left(  x_{0}\right)  }%
{\varphi\left(  x_{0}\right)  }%
\]
By the Harnack inequality, we have%
\[
u\left(  x\right)  \geq c\left(  \varepsilon+\frac{Md\left(  x_{0}\right)
}{\varphi\left(  x_{0}\right)  }\right)
\]
in $B_{\frac{d\left(  x_{0}\right)  }{2}}\left(  x_{0}\right)  $. Let $v$ be
the harmonic function in $B_{d\left(  x_{0}\right)  }\left(  x_{0}\right)
\backslash B_{\frac{d\left(  x_{0}\right)  }{2}}\left(  x_{0}\right)  $ such
that $v=0$ on $\partial B_{d\left(  x_{0}\right)  }\left(  x_{0}\right)  $ and
$v=\frac{cMd\left(  x_{0}\right)  }{\varphi\left(  x_{0}\right)  }$ on
$\partial B_{\frac{d\left(  x_{0}\right)  }{2}}\left(  x_{0}\right)  $, i.e.,%
\[
v\left(  x\right)  =\left\{
\begin{array}
[c]{lcc}%
\frac{cMd\left(  x_{0}\right)  }{\varphi\left(  x_{0}\right)  \left(
2^{n-2}-1\right)  }\left(  \left(  \frac{\left\vert x-x_{0}\right\vert
}{d\left(  x_{0}\right)  }\right)  ^{2-n}-1\right)  & \text{if} & n\geq3,\\
-\frac{cMd\left(  x_{0}\right)  }{\varphi\left(  x_{0}\right)  \ln2}\ln
\frac{\left\vert x-x_{0}\right\vert }{d\left(  x_{0}\right)  } & \text{if} &
n=2.
\end{array}
\right.
\]
It is easy to verify that for $x\in B_{d\left(  x_{0}\right)  }\left(
x_{0}\right)  $,
\[
v\left(  x\right)  \geq\frac{cM}{\varphi\left(  x_{0}\right)  }x_{1}+o\left(
\left\vert x\right\vert \right)  .
\]
Hence we have, near the origin,
\[
u\left(  x\right)  \geq v\left(  x\right)  \geq\frac{cM}{\varphi\left(
x_{0}\right)  }x_{1}+o\left(  \left\vert x\right\vert \right)  .
\]
in $B_{d\left(  x_{0}\right)  }\left(  x_{0}\right)  \backslash B_{\frac
{d\left(  x_{0}\right)  }{2}}\left(  x_{0}\right)  $. And for $x\notin
B_{d\left(  x_{0}\right)  }\left(  x_{0}\right)  $, we have
\[
x_{1}\leq\frac{\left\vert x\right\vert ^{2}}{2d\left(  x_{0}\right)  },
\]
and $\left(  \ref{Equation harmonic function positivity}\right)  $ follows
from $u\left(  x\right)  \geq0$. Now similar arguments as
\cite{2001Athanasopoulos_Caffarelli_Kenig_Salsa} imply%
\begin{equation}
\left\vert P\left(  A_{t},B_{3r^{\ast}}\left(  x^{\ast}\right)  \right)
-P\left(  A,B_{3r^{\ast}}\left(  x^{\ast}\right)  \right)  \right\vert
\leq\frac{C}{\varphi\left(  x_{0}\right)  }\left\vert V_{t}\right\vert
\end{equation}
and%
\begin{equation}
\int_{Z_{t}}\left\vert \triangledown u\right\vert ^{2}-\int_{Z_{t}}\left\vert
\triangledown v_{t}\right\vert ^{2}\geq\frac{cM^{2}}{\varphi\left(
x_{0}\right)  ^{2}}\left\vert V_{\frac{t}{2}}\right\vert .
\end{equation}
Summarize, we have shown that for $t$ sufficiently small,
\[
\frac{cM^{2}}{\varphi\left(  x_{0}\right)  ^{2}}\left\vert V_{\frac{t}{2}%
}\right\vert \leq\frac{C}{\varphi\left(  x_{0}\right)  }\left\vert
V_{t}\right\vert +C\left\vert V_{t}\right\vert .
\]
Finally, since $0\in\partial A$, Lemma \ref{Lemma positive density} implies
\[
\left\vert V_{t}\right\vert \geq ct^{n}%
\]
for some constant $c>0$, which guarantees the existence of a sequence of
positive numbers $\left\{  t_{j}\right\}  _{j=1}^{\infty}$, $t_{j}%
\rightarrow0$, such that
\[
\left\vert V_{t_{j}}\right\vert \leq2^{2n}\left\vert V_{\frac{t_{j}}{2}%
}\right\vert \text{.}%
\]
Let $t=t_{j}$ be sufficiently small, we have%
\[
\frac{cM^{2}}{\varphi\left(  x_{0}\right)  ^{2}}\leq2^{2n}\left(  \frac
{C}{\varphi\left(  x_{0}\right)  }+C\right)  \text{,}%
\]
hence $M\leq C$ for some constant independent of $\varepsilon$. The conclusion
of Lemma follows by letting $\varepsilon\rightarrow0$.
\end{proof}

By the standard covering argument, we have

\begin{corollary}
\label{Corollary linear growth}Let $K\subset\subset\Omega,$ then for some
positive constant $C_{K}$,%
\[
u\left(  x\right)  \leq C_{K}\operatorname*{dist}\left(  x,\partial A\right)
\]
holds for any $x\in K$.
\end{corollary}

The sublinear growth of $u$ near the free boundary implies the local Lipschitz
continuity of $u$:

\begin{theorem}
$u$ is locally Lipschitz continuous in $\Omega$.
\end{theorem}

\begin{proof}
We only need to consider the Lipschitz continuity of $u$ in any ball
$B_{r_{0}}\left(  x_{0}\right)  $ such that $x_{0}\in\partial A$ and%
\[
B_{2r_{0}}\left(  x_{0}\right)  \subset\subset\Omega.
\]
Let $K=B_{2r_{0}}\left(  x_{0}\right)  $, from Corollary
\ref{Corollary linear growth}, we have, for any $x\in B_{2r_{0}}\left(
x_{0}\right)  \backslash A$,
\[
u\left(  x\right)  \leq c_{K}\operatorname*{dist}\left(  x,\partial A\right)
.
\]
Let $x_{1},x_{2}\in B_{r_{0}}\left(  x_{0}\right)  $. If $x_{1},x_{2}\in A$,
then we have
\[
\left\vert u\left(  x_{1}\right)  -u\left(  x_{2}\right)  \right\vert
=\left\vert 0-0\right\vert =0.
\]
If $x_{1}\in A$, $x_{2}\in B_{r_{0}}\left(  x_{0}\right)  \backslash A$, then%
\[
\left\vert u\left(  x_{1}\right)  -u\left(  x_{2}\right)  \right\vert
=u\left(  x_{2}\right)  \leq c_{K}\operatorname*{dist}\left(  x_{2},\partial
A\right)  \leq c_{K}\left\vert x_{1}-x_{2}\right\vert .
\]
Similarly,
\[
\left\vert u\left(  x_{1}\right)  -u\left(  x_{2}\right)  \right\vert \leq
c_{K}\left\vert x_{1}-x_{2}\right\vert
\]
if $x_{2}\in A$, $x_{1}\in B_{r_{0}}\left(  x_{0}\right)  \backslash A$. So we
only need to consider the remaining case $x_{1},x_{2}\in B_{r_{0}}\left(
x_{0}\right)  \backslash A$, if%
\[
\left\vert x_{1}-x_{2}\right\vert \geq\frac{1}{2}\min\left\{
\operatorname*{dist}\left(  x_{1},\partial A\right)  ,\operatorname*{dist}%
\left(  x_{2},\partial A\right)  \right\}  ,
\]
without loss of generality, we assume%
\[
\left\vert x_{1}-x_{2}\right\vert \geq\frac{1}{2}\operatorname*{dist}\left(
x_{1},\partial A\right)  ,
\]
then we have%
\[
\left\vert u\left(  x_{1}\right)  -u\left(  x_{2}\right)  \right\vert \leq
u\left(  x_{1}\right)  \leq c_{K}\operatorname*{dist}\left(  x_{1},\partial
A\right)  \leq2c_{K}\left\vert x_{1}-x_{2}\right\vert .
\]
On the other hand, if%
\[
\left\vert x_{1}-x_{2}\right\vert \leq\frac{1}{2}\min\left\{
\operatorname*{dist}\left(  x_{1},\partial A\right)  ,\operatorname*{dist}%
\left(  x_{2},\partial A\right)  \right\}  ,
\]
then let
\[
r_{1}=\operatorname*{dist}\left(  x_{1},\partial A\right)  \leq\left\vert
x_{1}-x_{0}\right\vert <r_{0},
\]
we have%
\[
x_{2}\in B_{\frac{1}{2}r_{1}}\left(  x_{1}\right)  \subset B_{r_{1}}\left(
x_{1}\right)  \subset B_{2r}\left(  x\right)  \backslash A,
\]
since $u$ is harmonic in $B_{r_{1}}\left(  x_{1}\right)  \backslash A$, we
have%
\[
\left\Vert \triangledown u\right\Vert _{L^{\infty}\left(  B_{\frac{1}{2}r_{1}%
}\left(  x_{1}\right)  \right)  }\leq\frac{c\left(  n\right)  }{r_{1}%
}\left\Vert u\right\Vert _{L^{\infty}\left(  B_{\frac{r_{1}}{2}}\left(
x_{1}\right)  \right)  }\leq\frac{c\left(  n\right)  }{2}c_{K},
\]
here we have applied Corollary \ref{Corollary linear growth} in the last
inequality. Combining all the possibility, we have%
\[
\left\vert u\left(  x_{1}\right)  -u\left(  x_{2}\right)  \right\vert \leq
c_{1}\left(  n\right)  c_{K}%
\]
for any $x_{1},x_{2}\in B_{r_{0}}\left(  x_{0}\right)  $, hence $u$ is locally Lipschitz.
\end{proof}

\section{Analyticity of the reduced boundary}

The Lipschitz continuity of $u$ implies that $\partial A\cap\Omega$ is almost
area-minimizing in the sense introduced by F. J. Almgren \cite{1976Almgren}.

More precisely, we have

\begin{lemma}
For any $K\subset\subset\Omega$, there exists $C>0$, such that%
\begin{equation}
P\left(  A,B_{r}\left(  x\right)  \right)  \leq P\left(  \tilde{A}%
,B_{r}\left(  x\right)  \right)  +Cr^{n}
\label{Equation almost area-minimizing}%
\end{equation}
holds for every $x\in K$, every $r<\frac{1}{3}\operatorname*{dist}\left(
K,\partial\Omega\right)  $, and every $\tilde{A}$ which agrees with $A$ away
from $B_{r}\left(  x\right)  $.
\end{lemma}

\begin{proof}
Given $K\subset\subset\Omega$, let $x_{0}\in K$, $r_{0}<\frac{1}%
{3}\operatorname*{dist}\left(  K,\partial\Omega\right)  $ and $\tilde{A}$ be
any set which agrees with $A$ away from $B_{r_{0}}\left(  x_{0}\right)  $. We
can assume
\[
P\left(  A,B_{r_{0}}\left(  x_{0}\right)  \right)  >0,
\]
or else $\left(  \ref{Equation almost area-minimizing}\right)  $ holds
trivially. Since $u$ is Lipschitz continuous and $u\left(  x\right)  =0$ for
$x\in\partial A$, we have%
\[
u\left(  x\right)  \leq Cr_{0}%
\]
for any $x\in B_{2r_{0}}\left(  x_{0}\right)  $. Now let $\varphi\in
C^{\infty}\left(  \mathbb{R}^{n}\right)  $ be a function such that
$\varphi\equiv0$ in $B_{r_{0}}\left(  x_{0}\right)  $, $\varphi\equiv1$ away
from $B_{2r_{0}}\left(  x_{0}\right)  $ and%
\[
\left\Vert \varphi\right\Vert _{L^{\infty}}=1,\left\Vert \triangledown
\varphi\right\Vert _{L^{\infty}}\leq\frac{c\left(  n\right)  }{r_{0}}.
\]
Let%
\[
\tilde{u}=\varphi u,
\]
then
\[
\tilde{u}\left(  x\right)  =0\text{ for any }x\in\tilde{A}.
\]
Hence, the minimality of $\left(  A,u\right)  $ implies%
\[
\int_{B_{2r_{0}}}\left\vert \triangledown u\right\vert ^{2}+P\left(
A,B_{r_{0}}\left(  x\right)  \right)  \leq\int_{B_{2r_{0}}}\left\vert
\triangledown\tilde{u}\right\vert ^{2}+P\left(  \tilde{A},B_{r_{0}}\left(
x\right)  \right)  +C\left\vert \left\vert \tilde{A}\right\vert -\left\vert
A\right\vert \right\vert .
\]
Now%
\begin{align*}
&  \int_{B_{2r_{0}}}\left\vert \triangledown\tilde{u}\right\vert ^{2}%
=\int_{B_{2r_{0}}}\left\vert \triangledown\varphi u+\varphi\triangledown
u\right\vert ^{2}\\
\leq &  \left\vert B_{2r_{0}}\right\vert \left(  \left\Vert u\right\Vert
_{L^{\infty}\left(  B_{2r_{0}}\right)  }^{2}\cdot\left\Vert \triangledown
\varphi\right\Vert _{L^{\infty}\left(  B_{2r_{0}}\right)  }^{2}+\left\Vert
\varphi\right\Vert _{L^{\infty}\left(  B_{2r_{0}}\right)  }^{2}\cdot\left\Vert
\triangledown u\right\Vert _{L^{\infty}\left(  B_{2r_{0}}\right)  }^{2}\right)
\\
\leq &  Cr_{0}^{n},
\end{align*}
and%
\[
C\left(  \left\vert \tilde{A}\right\vert -\left\vert A\right\vert \right)
_{+}\leq Cr_{0}^{n},
\]
hence, we have%
\[
P\left(  A,B_{r_{0}}\left(  x\right)  \right)  \leq P\left(  \tilde
{A},B_{r_{0}}\left(  x\right)  \right)  +Cr_{0}^{n}.
\]
\end{proof}

Now the regularity result on almost area-minimizing boundaries implies
\cite{1976Almgren} \cite{1982Bombieri}:

\begin{theorem}
The reduced boundary $\partial^{\ast}A\cap\Omega$ is a $C^{1,\frac{1}{2}}$
hypersurface and the singular set $\left(  \partial A\backslash\partial^{\ast
}A\right)  \cap\Omega$ has Hausdorff dimension at most $n-8$.
\end{theorem}

To obtain higher order regularity of regular part of the free boundary, we
consider its Euler-Lagrange equation. Since $\partial^{\ast}A\cap\Omega$ is a
$C^{1,\frac{1}{2}}$ hypersurface, elliptic regularity theory implies that
$\triangledown u\in C^{\frac{1}{2}}\left(  \Omega\backslash\overline
{A}\right)  $ and we use $\triangledown u^{+}$ to denote its trace on
$\partial^{\ast}A\cap\Omega$, then $\triangledown u^{+}$ is $C^{\frac{1}{2}}$
on $\partial^{\ast}A\cap\Omega$.

\begin{lemma}
Free boundary equation%
\[
\kappa=\left\vert \triangledown u^{+}\right\vert ^{2}+C
\]
is satisfied weakly along $\partial^{\ast}A\cap\Omega$, where $C$ is constant
in any component of $\partial^{\ast}A\cap\Omega$.
\end{lemma}

\begin{proof}
We use the technique of domain variation. Without loss of generality, we
assume
\[
0\in\partial^{\ast}A\cap\Omega\text{, }B_{r}\left(  0\right)  \subset\Omega,
\]
and $B_{r}\left(  0\right)  \cap\partial A$ is the graph of $C^{1,\frac{1}{2}%
}$ function defined on a hyperplane passing through the origin. Let
$\varphi\in C_{0}^{\infty}\left(  B_{r}^{n-1}\left(  0\right)  ,\mathbb{R}%
^{n}\right)  $, we consider the mapping%
\[
T_{\varepsilon}:B_{r}\left(  0\right)  \rightarrow\mathbb{R}^{n}%
\]
given by
\[
y=T_{\varepsilon}\left(  x\right)  =x+\varepsilon\varphi.
\]
Let $\varepsilon<\operatorname*{dist}\left(  K,\partial B_{r}\left(  0\right)
\right)  $ be sufficiently small, where $K$ is the support of $\varphi$, then%
\[
T_{\varepsilon}:B_{r}\left(  0\right)  \rightarrow B_{r}\left(  0\right)
\]
is a diffeomorphism. We define $\left(  A_{\varepsilon},u_{\varepsilon
}\right)  $ such that
\[
A_{\varepsilon}\cap B_{r}\left(  0\right)  =T_{\varepsilon}\left(  A\cap
B_{r}\left(  0\right)  \right)  ,
\]
and%
\[
u_{\varepsilon}\left(  x\right)  =u\left(  T_{\varepsilon}^{-1}x\right)
\text{ for any }x\in B_{r}\left(  0\right)  .
\]
So from the minimality of $\left(  A,u\right)  $, we have%
\begin{align*}
&\int_{B_{r}}\left\vert \triangledown u\right\vert
^{2}dx+H^{n-1}\left(
\partial A\cap B_{r}\right)\\  \leq & \int_{B_{r}}\left\vert \triangledown
u_{\varepsilon}\right\vert ^{2}dx+H^{n-1}\left(  \partial A_{\varepsilon}\cap
B_{r}\right)  +C\left\vert \left\vert A_{\varepsilon}\right\vert -\left\vert
A\right\vert \right\vert ,
\end{align*}
where the last term comes from volume constraint.

First, we have%
\begin{align*}
&  \lim_{\varepsilon\rightarrow0}\frac{\int_{B_{r}}\left\vert \triangledown
u_{\varepsilon}\right\vert ^{2}dx-\int_{B_{r}}\left\vert \triangledown
u\right\vert ^{2}dx}{\varepsilon}\\
=  &  \int_{B_{r}}\left\vert \triangledown u\right\vert ^{2}%
\operatorname*{div}\varphi-2\triangledown u\triangledown\varphi\left(
\triangledown u\right)  ^{\intercal}\\
=  &  \int_{B_{r}\backslash\overline{A}}\left\vert \triangledown u\right\vert
^{2}\operatorname*{div}\varphi-2\triangledown u\triangledown\varphi\left(
\triangledown u\right)  ^{\intercal}.
\end{align*}
Integration by parts, using the fact that $\triangle u=0$ in $B_{r}\left(
0\right)  \backslash A$ and $u=0$ on $\partial A$, we have%
\begin{align*}
&  \int_{B_{r}\backslash\overline{A}}\left\vert \triangledown u\right\vert
^{2}\operatorname*{div}\varphi-2\triangledown u\triangledown\varphi\left(
\triangledown u\right)  ^{\intercal}\\
=  &  \int_{\partial A}\left\vert \triangledown u^{+}\right\vert ^{2}%
\varphi\cdot\nu-2\left(  \varphi\cdot\triangledown u^{+}\right)  \left(
\nu\cdot\triangledown u^{+}\right) \\
=  &  -\int_{\partial A}\varphi\cdot\nu\left\vert \triangledown u^{+}%
\right\vert ^{2}.
\end{align*}
Next,%
\[
\lim_{\varepsilon\rightarrow0}\frac{P_{\Omega}\left(  A_{\varepsilon}\right)
-P_{\Omega}\left(  A\right)  }{\varepsilon}=\int_{\partial A}\kappa
\varphi\cdot\nu,
\]
here $\int_{\partial A}\kappa\varphi\cdot\nu$ is well defined in weak sense
because the divergence structure of mean curvature, a precise formulation can
be given using local coordinates. Finally, we have%
\[
\lim_{\varepsilon\rightarrow0}\frac{\left\vert A_{\varepsilon}\right\vert
-\left\vert A\right\vert }{\varepsilon}=\int_{\partial A}\varphi\cdot\nu.
\]
Hence, we deduce%
\[
-\int_{\partial A}\varphi\cdot\nu\left\vert \triangledown u^{+}\right\vert
^{2}+\int_{\partial A}\kappa\varphi\cdot\nu+C\left\vert \int_{\partial
A}\varphi\cdot\nu\right\vert \geq0.
\]
If we change $\varphi$ to $-\varphi$, we have%
\[
-\int_{\partial A}\varphi\cdot\nu\left\vert \triangledown u^{+}\right\vert
^{2}+\int_{\partial A}\kappa\varphi\cdot\nu-C\left\vert \int_{\partial
A}\varphi\cdot\nu\right\vert \leq0.
\]
Hence for any $\varphi\in C_{0}^{\infty}\left(  B_{r}\left(  0\right)
,\mathbb{R}^{n}\right)  $ such that
\[
\int_{\partial A}\varphi\cdot\nu=0,
\]
we have%
\[
-\int_{\partial A}\varphi\cdot\nu\left\vert \triangledown u^{+}\right\vert
^{2}+\int_{\partial A}\kappa\varphi\cdot\nu=0.
\]
And we deduce
\[
-\left\vert \triangledown u^{+}\right\vert ^{2}+\kappa=C
\]
is satisfied weakly on $\partial A\cap B_{r}$, here%
\[
C=-\int_{\partial A}\varphi^{\ast}\cdot\nu\left\vert \triangledown
u^{+}\right\vert ^{2}+\int_{\partial A}\kappa\varphi^{\ast}\cdot\nu
\]
for any $\varphi^{\ast}$ satisfying%
\[
\int_{\partial A}\varphi^{\ast}\cdot\nu=0.
\]
\end{proof}

Now standard elliptic regularity theory implies $\partial^{\ast}A$ is a
$C^{2,\frac{1}{2}}$ hypersurface and the analyticity of $\partial^{\ast}A$
follows from a standard bootstrapping argument.

\section*{Acknowledgements}

The author would like to thank Professor Fanghua Lin for his constant help and encouragement.


\end{document}